\documentclass[12pt]{article}

\usepackage{amsmath,amsthm,amssymb}
\usepackage{graphics,graphicx}
\usepackage{e-jc}
\usepackage{color}
\usepackage{xspace}
\usepackage{pstricks,pst-node}
\usepackage[hypertex,colorlinks=true,linkcolor=green,filecolor=brown,citecolor=green]{hyperref}

\def\gl{ground level\xspace}
\def\gf{generating function\xspace}
\def\T{\ensuremath{\mathsf{T}}\xspace}
\def\U{\ensuremath{\mathsf{U}}\xspace}
\def\P{\ensuremath{\mathsf{P}}\xspace}

\def\mbf#1{\mathchoice{\hbox{\boldmath $\displaystyle #1$}}
        {\hbox{\boldmath $\textstyle #1$}}
        {\hbox{\boldmath $\scriptstyle #1$}}
        {\hbox{\boldmath $\scriptscriptstyle #1$}}} % mathboldface

\newtheorem{theo}{Theorem}

\title{Circular Digraph Walks, $k$-Balanced Strings, Lattice Paths and  Chebychev Polynomials}

\author{Evangelos Georgiadis\thanks{Corresponding Author}\\ \small Massachusetts Institute of Technology \\[-0.8ex] \small Cambridge, MA 02139, U.S.A. and\\[-.8ex] \small Dept. of Computer Science \& Tech, Tsinghua University\\[-.8ex] \small Beijing, 100084, P. R. China \\[-.8ex] \small \texttt{egeorg@mit.edu}\and David Callan\\ \small Department of Statistics\\[-.8ex] \small  University of Wisconsin-Madison\\[-0.8ex] \small 1300 University Ave\\[-.8ex]\small Madison, WI 53706-1532, U.S.A.\\[-0.8ex]\small \texttt{callan@stat.wisc.edu}\and Qing-Hu Hou\thanks{This work was supported by the 973 Project, the PCSIRT Project of the Ministry of Education, the Ministry of Science and Technology,
and the National Science Foundation of China.}\\\small Center for Combinatorics, LPMC-TJKLC\\[-.8ex] \small Nankai University\\[-.8ex] \small Tianjin 300071, P. R. China \\[-.8ex] \small \texttt{hou@nankai.edu.cn}}

\date{\dateline{June 26, 2008}{XXX}}

\begin{document}
  \maketitle  
\begin{abstract}
We count the number of walks of length $n$ on a $k$-node circular 
digraph 
that cover all $k$ nodes in two ways. The first way illustrates the transfer-matrix 
method. The second involves counting various classes of height-restricted lattice paths.  
We observe that the results also count so-called $k$-balanced strings of length
$n$, generalizing a 1996 Putnam problem.

\end{abstract}

\section{Introduction: Walks and $\mbf{k}$-Balanced Binary Strings}

\begin{figure}[ht]
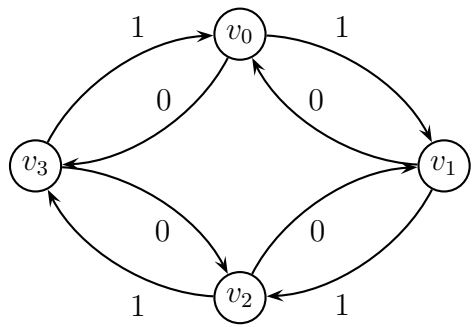
\label{demonstration1}
\centering $\psmatrix[colsep=2cm,rowsep=1cm,mnode=circle]
&v_0\\
v_3&&v_1\\
&v_2
\psset{arrowscale=1.5} 
\ncarc[arcangle=30]{->}{1,2}{2,1}<{0}
\ncarc[arcangle=-30]{<-}{1,2}{2,1}^{1}
\ncarc[arcangle=30]{->}{1,2}{2,3}^{1}
\ncarc[arcangle=-30]{<-}{1,2}{2,3}>{0}
\ncarc[arcangle=30]{->}{2,3}{3,2}_{1}
\ncarc[arcangle=-30]{<-}{2,3}{3,2}>{0}
\ncarc[arcangle=30]{->}{2,1}{3,2}<{0}
\ncarc[arcangle=-30]{<-}{2,1}{3,2}_{1}
\endpsmatrix
$ \caption{\small When $k=4$, an instance of a good walk of length
$5$ starting from $v_0$ is $(v_0,v_1,v_2,v_1,v_2,v_3)$. This walk
generates the unique binary string $11011$.}
\end{figure}

Let $C_k$ be a circular digraph that consists of $k$ nodes, namely,
$v_{0},\ldots, v_{k-1}$. A walk on $C_k$ of length $n$ is simply a
sequence of $n+1$ nodes $(w_0,\ldots,w_n)$ such that $w_i$ is
adjacent to $w_{i+1}$ in $C_k$ for $0 \le i \le n-1$. Notice that we
may assign the (clockwise) arcs, between nodes $v_i$ and
$v_{{(i+1)}{\pmod{k}}}$ for each $i=0, \ldots, k-1$, with transition
label $1$ whereas assign the (counterclockwise) arcs, between
$v_{{(i+1)}{\pmod{k}}}$ and $v_{i}$ for each $i=0, \ldots, k-1$,
with transition label $0$. Then each walk on $C_k$ of length $n$
generates a unique binary word of length $n$. For ease of
visualization, we provide Figure \ref{demonstration1} as an instance
when $k=4$.

We now define a ``good walk'' on $C_k$ as a walk starting from $v_0$
and visiting all $k$ nodes of $C_k$. We settle the question of how
many good walks exist, by restricting our attention to ``bad walks''
(i.e. walks that do not cover all nodes).

The binary strings generated by ``bad walks'' can be placed into a
$1$-$1$ correspondence with the so-called $(k-2)$-balanced binary
strings. A $k$-balanced binary string, in turn, is defined as a
finite binary string $S$ in which every substring $T$ (of
consecutive bits) of $S$ has $-k \leq \Delta(S) \leq k$, where
$\Delta(S)$ denotes the number of $1$'s minus the number of $0$'s.
For example, $11011$ represents an unbalanced string (for $2$-balanced binary strings).

A $1996$ Putnam problem \cite{KPV} by Michael Larsen asked 
for the number of $2$-balanced binary strings, and a generalization 
to $k$-balanced strings was the motivation for this paper. An 
explicit-sum solution to the Putnam problem is given in \cite{KPV} but 
generalizing it seems unwieldy. Here we focus on generating 
functions. 

The outline of the paper is as follows. In Section 2 we 
use the transfer-matrix method to  obtain the desired \gf as a difference 
$S_{k+1}(x)-S_{k}(x)$ where $S_{k}(x)$ is the sum of the entries in a 
certain $k\times k$ matrix, and to make a first stab at simplifying 
$S_{k}(x)$. In Section 3 we survey the use of Chebychev polynomials 
to count various classes of height-restricted lattice paths and deduce 
an alternative expression for the desired \gf as a product 
$R_{k}(x)R_{k-1}(x)$. In Section 4 we reconcile the two formulas 
$S_{k+1}(x)-S_{k}(x)$ and $R_{k}(x)R_{k-1}(x)$.

\section{The Transfer-Matrix Approach}

\subsection{The basic result}\label{parta}

\begin{theo}
Let $A_k$ denote the tridiagonal $k \times k$
matrix with 1s just above and just below the main diagonal and 0s 
elsewhere,
\begin{equation}
\label{tri-diag}
A_k = \left( \begin{array}{ccccc} 0 & 1 & 0 & \cdots & 0 \\
1 & 0 & 1 & \cdots & 0 \\
0 & 1 & 0 & \ddots & \vdots  \\
\vdots & \ddots & \ddots & \ddots & 1 \\
0 & \cdots & 0 & 1 & 0
\end{array} \right),
\end{equation}
and let $S_k(x)$ denote the sum of all the entries in $(I_{k}-xA_{k})^{-1}$ 
where $I_{k}$ is the $k \times k$ identity matrix. 
Then the generating function
for ``bad walks'' of length $n$ on $C$ equals
\[
S_{k-1}(x)-S_{k-2}(x).
\]
In other words, the generating function for $k$-balanced strings of
length $n$ is
\[
f_k(x) = S_{k+1}(x) - S_k(x).
\]
\end{theo}

\noindent \textbf{Proof}\quad Given a ``bad walk'' $w=(w_0=v_0,w_1,\ldots,w_n)$ of length $n$,
let
\[
\max (w) = \max \{i \colon \{v_0,v_1,\ldots,v_i\} \subseteq \{w_0,
\ldots, w_n \} \}.
\]
We see that ``bad walks'' with $\max(w)=r$ are just the walks
$w=(w_0,\ldots,w_n)$ on $C \setminus \{v_{r+1}\}$ such that
$w_0=v_0$ and $v_r \in \{w_0,\ldots,w_n \}$. Notice that an
arbitrary walk $w=(w_0,\ldots,w_n)$ on $C \setminus \{v_{r+1}\}$
either satisfies $v_r \in \{w_0,\ldots,w_n \}$ or is a walk on $C
\setminus \{v_r, v_{r+1}\}$. Thus the walks with $\max(w)=r$ are those 
that miss $\{v_{r+1}\}$ but don't miss $\{v_r, v_{r+1}\}$.
Now $C\setminus \{v_{r+1}\}$ is the path graph on vertex list 
$v_{r+2},v_{r+3},\ldots,v_{k-1},v_{0},v_{1},\ldots,v_{r}$. Note that 
$v_{0}$ is the $(k-1-r)$th entry in the vertex list and the adjacency 
matrix is $A_{k-1}$. The transfer-matrix method
\cite[Theorem 4.7.2]{Stanley} says that the generating
function for walks from $v_{0}$ to the $j$th vertex is the $(k-1-r,j)$ 
entry of $(I_{k-1}-xA_{k-1})^{-1}$. Similarly, the generating
function for walks from $v_{0}$ to the $j$th vertex in the path graph $C
\setminus \{v_r, v_{r+1}\}$ is the $(k-1-r,j)$ 
entry of $(I_{k-2}-xA_{k-2})^{-1}$. Taking the difference and summing over $r$ and $j$ yields the 
result. \qed

\subsection{A determinant formula}\label{partb}

Now we obtain an expression for $S_k(x):=$ sum of entries in 
$(I_{k}-xA_{k})^{-1}$. Let $\U_{k}(x)$ denote the 
``combinatorial'' Chebyshev polynomial introduced in the next section.

Since $\det(I_{k}-xA_{k})$ and $\U_{k}(x)$ both satisfy the recurrence
$\P_{k}(x)=\P_{k-1}(x)-x^{2}\P_{k-2}(x)$ with initial conditions 
$\P_{0}(x)=\P_{1}(x)=1$, we conclude that $\det(I_{k}-xA_{k})=\U_{k}(x)$.

Applying linear algebra, we have $S_k(x) =  x_1+\cdots+x_k$, 
where the $x_i$ (functions of $x$) denote the solutions to the equation system
\[
(I_{k}-xA_{k}) \left(\begin{array}{c} x_1 \\ x_2 \\ \vdots \\x_k \end{array}
\right) = \left(\begin{array}{c} 1 \\ 1 \\ \vdots \\1
\end{array} \right).
\]
By summing up the equations, we find
\begin{equation}
    (1-x)x_1 + (1-2x)(x_2 + \cdots + x_{k-1}) + (1-x)x_k = k,
    \label{eq:1}
\end{equation}
and Cramer's rule implies
\begin{equation}
    x_1 = x_k = \left. \det \left( \begin{array}{cccccc} 1 & -x & 0 & 0 & \cdots & 0 \\
1 & 1 & -x & 0 & \cdots & 0 \\
1 & -x & 1 & -x & \cdots & \vdots  \\
1 & 0 & -x & 1 & \ddots &  \vdots  \\
\vdots & \vdots & \ddots & \ddots & \ddots& \vdots \\
1 & 0 & \cdots & 0 & -x & 1
\end{array} \right) \right/ \U_k(x).
    \label{eq:2}
\end{equation}
Denote the determinant in the numerator by $W_k$. Thus from 
(\ref{eq:1}) and (\ref{eq:2}), we have
\[ 
S_k(x) = \frac{k - 2 x x_1}{1 - 2 x}
=\frac{k \U_k(x) - 2 x W_k}{(1 - 2 x)\U_{k}(x)}.
\]

\vspace{2mm}

\section{The Lattice Path Approach}

\subsection{Combinatorial Chebychev polynomials}

The familiar Chebychev polynomials $T_{k}(x)$ 
and $U_{k}(x)$ (first and second kinds) are defined by 
$\cos k \theta = T_{k}(\cos \theta)$ and $\sin 
(k+1) \theta /\sin \theta = U_{k}(\cos \theta)$. They occur in diverse 
areas, as suggested by the subtitle of Theodore Rivlin's book 
\cite{rivlin90}. Their application in combinatorics to lattice path counting is less 
well known. For this purpose, it is convenient 
to define modified Chebychev polynomials by
\[
\T_{k}(x) = 2x^{k}T_{k}\big(\textrm{\footnotesize{$\frac{1}{2x}$}} \big), \qquad 
\U_{k}(x) =x^{k}U_{k}\big(\textrm{\footnotesize{$\frac{1}{2x}$}} \big).
\]
This removes an extraneous power of 2 and 
reverses the coefficients to produce integer-coefficient  
polynomials with constant term 1 (except that $\T_{0}=2$) which might be called the 
\emph{combinatorial Chebychev polynomials}. Both satisfy the defining 
recurrence $\P_{k}(x)=\P_{k-1}(x)-x^{2}\P_{k-2}(x)$, differing only in 
the initial conditions, and both have simple explicit 
expressions:
\[
\T_{k}(x)=\sum_{j=0}^{\lfloor k/2 
\rfloor}(-1)^{j}\left(\binom{k-j}{j}+\binom{k-j-1}{j-1}\right)x^{2j}, \quad 
\U_{k}(x)=\sum_{j=0}^{\lfloor k/2 \rfloor}(-1)^{j}\binom{k-j}{j}x^{2j}.
\]
The first few are listed in the following Table.

\vspace*{-6mm}

\[
\begin{array}{cccc}
    k & & \T_{k}(x)  &\U_{k}(x) \\ \hline
    0 &  & 2  & 1\\
    1 &  & 1  & 1\\ 
    2 &  & 1 - 2x^{2}  & 1-x^{2} \\ 
    3 & &  1 - 3x^{2} & 1 - 2x^{2}\\ 
    4 & &  1 - 4x^{2} + 2x^4  &  1 - 3x^{2} + x^4\\
    5 & &  1 - 5x^{2} + 5x^4 & 1 - 4x^{2} + 3x^4\\
    6 & &  1 - 6x^{2} + 9x^4 - 2x^6 &  1 - 5x^{2} + 6x^4 - x^6\\
    7 & &  1 - 7x^{2} + 14x^4 - 7x^6 & 1 - 6x^{2} + 10x^4 - 4x^6

\end{array}
\]

\vspace*{2mm}

\centerline{{\small Table of combinatorial Chebychev polynomials}}

\vspace*{4mm}

\subsection{Application to height-restricted lattice paths}

Consider lattice paths of upsteps $u=(1,1)$ and downsteps $d=(1,-1)$.
The horizontal line through a path's initial vertex is \emph{ground level} and 
heights are measured relative to \gl. Thus if the path starts at 
the $x$-$y$ origin, ground level is the $x$-axis. The \emph{height} of 
a path is the maximum of the heights of its vertices. A \emph{nonnegative} 
path is one that never dips below \gl. A \emph{balanced} path (not to 
be confused with $k$-balanced strings) is one 
that ends at \gl.
A \emph{Dyck} path is a nonnegative balanced path, 
including the empty path. 

The generating function for a given class of paths is $\sum_{n\ge 
0}a(n)x^{n}$ where $a(n)$ is the number of paths of size $n$: size is 
taken as ``number of steps'' except for paths specified to 
terminate at height $k$, where size is ``$\# \textrm{\,steps} -k$'' 
since such a path necessarily contains $k$ upsteps.

The basic application of combinatorial Chebychev polynomials to count 
height-restricted lattice paths is given in Table 1. \emph{Here, and in the 
sequel, $\U_{k}$ is short for  $\U_k(x)$ and so on.}
\vspace*{6mm}

\noindent \qquad \begin{tabular}{cc}
\multicolumn{2}{c}{ \framebox[\width][c]{ paths bounded by  $y=0$ and 
$y=k$ } }  \\[3mm]
path ends  & generating  \\
 at height   & function   \\ \hline  
\rule[-0.6cm]{0mm}{1.5cm}
$0$  & $F_{k} =$ 
{\Large $\frac{\vphantom{\int}\U_{k}}
{\vphantom{\int}\U_{k+1}}$} \\
\rule[-0.6cm]{0mm}{1.3cm}
$k$  &  $G_{k} =$
{\Large $\frac{1}
{\vphantom{\int}\U_{k+1}}$} 
\end{tabular} 
\hspace*{15mm}
\begin{tabular}{cc}
\multicolumn{2}{c}{ \framebox[\width][c]{ paths bounded by $y=\pm k$ } } 
\\[3mm]
path ends  & generating  \\
at height  & function   \\ \hline
\rule[-0.6cm]{0mm}{1.5cm}
$0$  &  $\overline{F}_{k} =$
{\Large $\frac{\vphantom{\int}\U_{k}}
{\vphantom{\int}\T_{k+1}}$} \\
\rule[-0.6cm]{0mm}{1.3cm}
$k$ &  $\overline{G}_{k} =$
{\Large $\frac{1}
{\vphantom{\int}\T_{k+1}}$}
\end{tabular}

\vspace*{6mm}

\centerline{Table 1}
\centerline{Generating functions for some 
height-restricted}\vspace*{-1mm}
\centerline{lattice paths with specified terminal height}
\vspace*{4mm}

Thus the first item, $F_{k}(x)$, is the \gf for Dyck paths of height 
$\le k$ with $x$ marking length. The expressions for $F_{k}$ 
and $G_{k}$ are folklore; two early references are 
\cite{kreweras70,knuth72} and a recent one is \cite{bousquet08}. For 
completeness we briefly outline below proofs for all the items in 
Table 1. 

It is also possible to find corresponding generating functions 
$H_{k}(x)$ and $\overline{H}_{k}(x)$ for paths 
with no restriction on the height of the terminal vertex. 
%Let $H_{k}$ denote 
%the GF for paths bounded by $y=0$ and $y=k$, and $R_{k}$ the GF for paths bounded by $y=\pm k$.
Here it is necessary to distinguish the cases $k=2m$ is even and 
$k=2m+1$ is odd:

\vspace*{6mm}

\noindent \quad  \begin{tabular}{c}
{ \framebox[\width][c]{ paths bounded by  $y=0$ and $y=k$ } } \\[1mm]
\rule[-0.6cm]{0mm}{1.8cm}
 $H_{2m} =$ 
{\Large $\frac{\vphantom{\int}\U_{m} + x 
\U_{m-1} }
{\vphantom{\int}\T_{m+1}}$} \\
\rule[-0.8cm]{0mm}{2cm}
 $H_{2m+1} =$
{\Large $\frac{\vphantom{\int} 
\U_{m}}
{\vphantom{\int}\U_{m+1} - x
\U_{m}}
$} 
\end{tabular} \hspace*{12mm}
\begin{tabular}{c}

{ \framebox[\width][c]{ paths bounded by $y=\pm k$ } } \\[1mm]
\rule[-0.6cm]{0mm}{1.8cm}
$\overline{H}_{2m} =$
{\Large $\frac{\vphantom{\int}\left(\U_{m} + x 
\U_{m-1}\right)^{2} }
{\vphantom{\int}\T_{2m+1}}$} \\
\rule[-0.8cm]{0mm}{2cm}
$\overline{H}_{2m+1} = (1+2x)$
{\Large $\frac{\vphantom{\int}
\left(\U_{m}\right)^{2}}
{\vphantom{\int}\T_{2m+2}}$}  
\end{tabular}
\vspace*{3mm}

\centerline{Table 2}
\centerline{Generating functions for some height-restricted lattice 
paths}

\vspace*{6mm}
Sri Gopal Mohanty \cite{mohanty79} uses the reflection principle to 
count paths bounded by $y=s$ and $y=-t$ for 
arbitrary nonnegative $s$ and $t$, obtaining explicit sums rather 
than generating functions.

\noindent \textbf{Proofs for Tables 1 and 2}

$\mbf{F}$\,:\quad A nonempty Dyck path $P$ can be uniquely expressed as $uP_{1}dP_{2}$ 
where $P_{1}$ and $P_{2}$ are Dyck paths. The path $P$ has height $\le 
k$ if and only if $P_{1}$ has height $\le k-1$ and $P_{2}$ has 
height $\le k$. This observation translates to a recurrence for the 
\gf:
\[
F_{k}=1+x^{2}F_{k-1}F_{k},
\]
with solution $F_{k}=\U_{k}/\U_{k+1}$ because the 
substitution $\U_{k}/\U_{k+1}$ for $F_{k}$ reduces to 
$\U_{k+1}=\U_{k}-x^{2} \U_{k-1}$, equivalent to a well known 
recurrence for Chebychev polynomials.

$\mbf{G}$\,:\quad A path bounded by $y=0$ and $y=k$ that terminates at height $k$ has a last upstep to height $i$ for 
$i=1,2,\ldots,k-1$ and the last upstep to height $k$ is necessarily
the last step of the path. Remove these $k$ upsteps to obtain a list of 
$k$ Dyck paths (some may be empty). The $i$th path in this list from right to left has height $\le i$ 
and hence \gf $F_{i}$. Thus $G_{k} = 
\prod_{i=1}^{k}F_{i}=1/\U_{k+1}$.

$\mbf{\overline{F}}$\,:\quad A balanced path $P$ bounded by $y=\pm k$ is either (i) empty or (ii) 
starts up or (iii) starts down. In case (ii) P decomposes as 
$uP_{1}dP_{2}$ where $P_{1}$ is a Dyck path of height $\le k-1$ and 
$P_{2}$ is another balanced path  bounded by $y=\pm k$. Thus case (ii) 
contributes $x^{2}F_{k-1}\overline{F}_{k}$ and, by symmetry, so does case 
(iii). Hence 
\[
\overline{F}_{k}=1+2x^{2}F_{k-1}\overline{F}_{k}
\]
leading to $\overline{F}_{k}=\U_{k}/(\U_{k}-2x^{2}\U_{k+1})$ 
and so, using another well known Chebychev 
polynomial identity, $\overline{F}_{k}=\U_{k}/\T_{k+1}$.

$\mbf{\overline{G}}$\,:\quad A path bounded by $y=\pm k$ terminating at height $k$ has a last 
upstep to height $i,\ 1\le i\le k-1$. Delete these upsteps to obtain a list 
consisting of a balanced path bounded by $y=\pm k$, followed by $k-1$ Dyck 
paths of heights $\le k-1,\ \le  k-2,\ \ldots,\ \le 1$ respectively. 
Thus 
\[
\overline{G}_{k}=\frac{\U_{k}}{\T_{k+1}}\frac{\U_{k-1}}{\U_{k}} 
\cdots \frac{\U_{1}}{\U_{2}} =\frac{1}{\T_{k+1}}.
\]

$\mbf{H}$\,:\quad A path bounded by $y=0$ and $y=k$ with no restriction on the terminal 
height is either (i) empty or (ii) starts with an upstep and never 
returns to \gl or (iii) has the form $uPdQ$ where $P$ is a Dyck path 
of height $\le k-1$ and $Q$ is a path bounded by $y=0$ and $y=k$. The 
contributions to the \gf $H_{k}$ are respectively (i) 1, (ii) 
$xH_{k-1}$, (iii) $x^{2}F_{k-1}H_{k}$. Thus 
\[
H_{k} = 1 + xH_{k-1}+x^{2}F_{k-1}H_{k}.
\]
It is routine, if tedious, to verify that the expressions for 
$H_{2m}$ and $H_{2m+1}$ in Table 2 satisfy this equation. The 
$2m$ case, for example, can be verified as follows. Replace 
$\T_{m+1}$ by $\U_{m+1}-x^{2} \U_{m-1}$ (another Chebychev 
identity) and revert to the 
standard Chebychev polynomials $U_{k}(x)$. With the substitution $y=1/(2x)$ this 
reduces matters to verifying that 
\[
\big(U_{m}^{2}(y)-U_{m-1}^{2}(y)\big) \big(2y 
U_{2m}(y)-U_{2m-1}(y)\big) = 
U_{m}(y)U_{2m}(y)\big(U_{m+1}(y)-U_{m-1}(y)\big),
\]
an identity that ultimately depends on the elementary addition formulae for 
trigonometric functions.

$\mbf{\overline{H}}$\,:\quad A path bounded by $y=\pm k$ with no restriction on the terminal 
height is either (i) empty or (ii) starts with an upstep (resp. 
downstep) and never returns to \gl or (iii) starts with an upstep (resp. 
downstep) and returns to \gl. Case (ii) ``start up'' makes a contribution of $x 
H_{k-1}$ and by symmetry, case (ii) ``start down'' makes the same 
contribution. In case (iii) ``start up'', the path has the form $uPdQ$  
where $P$ is a Dyck path of height $\le k-1$ and $Q$ is another path 
of the kind being counted. Thus case (iii) makes contribution $2x^{2} 
F_{k-1}\overline{H}_{k}$. Hence
\[
\overline{H}_{k}=1+2xH_{k-1}+2x^{2}F_{k-1}\overline{H}_{k}
\]
and another trite calculation shows that the expression for $\overline{H}_{k}$ 
in Table 2 satisfies this recurrence.

\subsection{Application to $\mbf{k}$-balanced strings}\label{R}
A binary string of, say, $X$s and $O$s can be coded as a lattice 
path: $X\rightarrow u,\: O\rightarrow d$. The $k$-balanced strings of 
length $n$ translate to lattice paths of $n$ steps with \emph{vertical 
extent} $\le k$ where vertical extent means 
``maximum vertex height $-$ minimum  vertex height''. A 
recurrence  for the \gf $g_{k}(x)$ for these paths (with $x$ marking 
number of steps) can be obtained from the following  decomposition. 
Such a path is either nonnegative or else dips below \gl and 
hence has a first downstep $d$ carrying it to its lowest level (below 
\gl). The first case gives contribution $H_{k}$. In the second 
case, the path has the form $PdQ$ where the reverse of $P$ is a 
nonnegative path of height $\le k-1$ and $Q$ is a nonnegative path of height 
$\le k$ as illustrated.
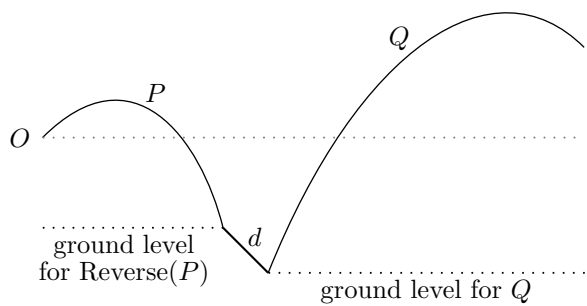
\begin{center} 
\begin{pspicture}(-4.5,-0.3)(12,3.5)
\psset{unit=.6cm}
    
\psbezier[linewidth=.5pt](0,3)(2,5)(3.5,3)(4,1)
\psbezier[linewidth=.5pt](5,0)(7,5)(10,7)(12,5)

\psline(4,1)(5,0)

\psline[linestyle=dotted](0,1)(4,1)
\psline[linestyle=dotted](5,0)(12,0)

\psline[linecolor=gray,linestyle=dotted](0,3)(12,3)

\rput(4.7,.8){\textrm{{\footnotesize $d$}}}
\rput(-.5,3){\textrm{{\footnotesize $O$}}}
\rput(2.5,4){\textrm{{\footnotesize $P$}}}
\rput(7.9,5.2){\textrm{{\footnotesize $Q$}}}
\rput(8.5,-0.4){\textrm{{\footnotesize ground level for $Q$}}}

\rput(1.8,0.6){\textrm{{\footnotesize ground level}}}
\rput(1.8,0.0){\textrm{{\footnotesize for Reverse($P$)}}}

\end{pspicture}
\end{center} 
Hence
\[
g_{k}=H_{k}+xH_{k-1}H_{k}.
\]
This is the desired \gf but it has an interesting alternative 
expression. 
Define a sequence of rational functions $(R_{k}(x))_{k\ge 1}$ by 
\[
R_{2m}=
\textrm{ $\frac{\vphantom{\int}\U_{m}}
{\vphantom{\int}\T_{m+1}}$},\quad R_{2m+1}=\textrm{  $\frac{\vphantom{\int} 
\vphantom{\int}\U_{m+1} + x \U_{m}}
{\vphantom{\int}\U_{m+1} - x \U_{m}} $}.
\]
Then it is easy to check that $H_{k}\big(1+xH_{k-1}\big) = 
R_{k}R_{k-1},\ k\ge 1$. Thus $g_{k}=R_{k}R_{k-1}$ and we have 
established
\begin{theo}
    The generating function for $k$-balanced binary strings, 
    equivalently for $u$-$d$ paths of vertical extent $\le k$, is 
    given by
\[
\begin{cases} \textrm{ {\large
  $ \frac{\vphantom{\int}\U_{m}}{\vphantom{\int}\T_{m+1}}\cdot
  \frac{\vphantom{\int}\U_{m} + x \U_{m-1} 
  }{\vphantom{\int}\U_{m}-x\U_{m-1}} $}} & 
   \text{if $k=2m$ is even;} \\[5mm]
   \textrm{ {\large
  $\frac{\vphantom{\int}\U_{m+1} + x \U_{m} 
  }{\vphantom{\int}\U_{m+1}-x\U_{m}}\cdot
   \frac{\vphantom{\int}\U_{m}}{\vphantom{\int}\T_{m+1}} $}} & \text{if $k=2m+1$ is odd.} 
\end{cases}   
\]
\end{theo}
\textbf{Remark}\quad The expression in Theorem 2 for $g_{k}=R_{k}R_{k-1}$ is 
in lowest terms because $\T_{m+1}=\U_m - 2x^2 \U_{m-1}$ and the 
recurrence $\U_{m}=\U_{m-1}-x^{2} \U_{m-2}$ yields by induction that $\U_{m}$ and $\U_{m-1}$ 
are relatively prime.

\textbf{Remark}\quad $R_{k}$ can be compactly expressed in terms of entries in 
Tables 1 and 2:
\[
R_{2m}=\overline{F}_{m},\quad R_{2m+1}=1+2xH_{2m+1}.
\]
Thus $g_{k}$ involves convolutions of paths bounded by $y=\pm 
\lfloor k/2 \rfloor$ terminating at \gl ($\overline{F}_{m}$) and nonnegative paths of 
height $\le k$ terminating anywhere ($H_{2m+1}$). 
A combinatorial explanation would be interesting but does not seem 
to be obvious.

\section{Reconciling the Two Formulas}
We have obtained expressions $f_{k}(x)$ and $g_{k}(x)$ for the \gf 
for $k$-balanced strings in Sections 2 and 3 respectively. We now 
show that $f_{k}=g_{k}$. The proof ultimately depends on the standard 
identities
\begin{equation}\label{basic}
    \begin{aligned}
    \U_{2k} &=  \U_{k}^{2}-x^{2}\U_{k-1}^{2},  \\
    \U_{2k+1} &=  \U_{k}^{2}-2x^{2}\U_{k}\U_{k-1}.
    \end{aligned}
\end{equation}
Set 
\[
P_k = \frac{k \U_k - 2 x W_k}{(1 - 2 x)}
\]
so that, from Section 
2.2, $S_{k}=P_{k}/\U_{k}$.

First, we find an expression for the determinant $W_{k}$. 
By cofactor expansion along the first row, $W_{k}$ satisfies the 
defining recurrence
\[
W_{0}=0,\ W_{1}=1,\quad W_{k}=\U_{k-1}+xW_{k-1},
\]
with solution, verified using (\ref{basic}),
\begin{eqnarray*}
    W_{2m} & = & (\U_{m} + x\U_{m-1})\U_{m-1},  \\
    W_{2m+1} & = & (\U_{m} + x\U_{m-1})\U_{m}.
\end{eqnarray*}
Now define two sequences of 
polynomials $(A_{k})_{k\ge 0},\ (B_{k})_{k\ge 0}$ by
\begin{eqnarray*}
    A_{2m} = \U_{m}-x\U_{m-1}, & \hspace*{5mm} & A_{2m+1}=\T_{m+1}= 
    \U_{m}-2x^{2}\U_{m-1}; \\
    B_{2m}=\U_{m}+x\U_{m-1}, & \hspace*{5mm} & B_{2m+1}=\U_{m}.
\end{eqnarray*}
Thus, in particular, $W_{k}=B_{k}B_{k-1}$ whether $k$ is even or odd.
Also define a sequence $(C_{k})_{k\ge 0}$ of rational functions 
(actually polynomials) by
\[
C_{k}=\frac{kA_{k}-2xB_{k-1}}{1-2x}.
\]
It is now easy to verify that 
\begin{eqnarray*}
    P_{k} & = & B_{k}C_{k},  \\
    \U_{k} & = & A_{k} B_{k}, \\
    R_{k} & = & \frac{B_{k+1}}{A_{k+1}},
\end{eqnarray*}
where $R_{k}$ is as defined in the preceding section, and that
\begin{equation*}
 C_{k}A_{k+1}-C_{k-1}A_{k}\  =  \  W_{k}. 
\end{equation*}
Hence
\[
f_{k}=S_{k+1}-S_{k}=\frac{P_{k+1}}{\U_{k+1}}-\frac{P_{k}}{\U_{k}}=
\frac{C_{k+1}}{A_{k+1}}-\frac{C_{k}}{A_{k}} = 
\frac{W_{k+1}}{A_{k+1}A_{k}} = 
\frac{B_{k+1}B_{k}}{A_{k+1}A_{k}}=R_{k}R_{k-1}=g_{k},
\]
as required.

\paragraph*{Acknowledgment}
One of us (E.G.) would like to thank Mike Sipser of MIT for great support, encouragement and inspiration 
as well as John Tsitsiklis of MIT. Two of us (E.G., Q-H.H.) would like to thank Dr. G.C. Xin for valuable comments.

\end{document}